\documentclass[a4paper,11pt]{amsart}
\usepackage{color}

\newcommand{\aspas}[1]{``{#1}''}
\usepackage{amsmath,amsfonts,amsthm,amstext}
\usepackage{amsmath,amssymb,indentfirst,epsfig,multicol}
\usepackage{amssymb,latexsym}
\usepackage{lipsum}
\usepackage{epsfig}
\usepackage{graphicx}
\usepackage{amscd}
\usepackage{verbatim}
\usepackage[normalem]{ulem}
\usepackage[english]{babel}
\usepackage[latin1]{inputenc}
\usepackage[all]{xy}

\newcommand{\cqd}{\hfill $\rule{2.5mm}{2.5mm}$}

\newtheorem{theorem}{Theorem}[section]
\newtheorem{rema}[theorem]{Remark}
\newtheorem{Que}[theorem]{Question}
\newtheorem{lema}[theorem]{Lemma}
\newtheorem{exam}[theorem]{Example}
\newtheorem{defi}[theorem]{Definition}
\newtheorem{prop}[theorem]{Proposition}
\newtheorem{cor}[theorem]{Corollary}

\def\proof{\vskip .2cm\noindent{\sc Proof}\quad}

\begin{document}

\title[On the topology of the Milnor fibration]{On the topology of the Milnor fibration}

\author{Souza, T. O.}
\address{Faculty of Mathematics, Federal University of Uberl\^andia, Uberl\^andia, Minas Gerais, Brazil.}
\email{tacioli@ufu.br}

\author{Zapata, C. A. I.}
\address{Departamento de Matem\'{a}tica, Universidade de S\~{a}o Paulo, Instituto de Matem\'{a}tica e  Estad\'{i}stica-IME/USP, R. do Mat\~{a}o, 1010 - Butant\~{a}, CEP:
05508-090 - S\~{a}o Paulo, Brasil}
\email{cesarzapata@usp.br}

\date{}



\thanks {The second author would like to thank grants \#2016/18714-8, \#2018/23678-6 and \#2022/03270-8, S\~{a}o Paulo Research Foundation (FAPESP) for financial support.}

\maketitle

\begin{abstract}
In this paper we present new results about the topology of the Milnor fibrations of analytic function-germs with a special attention to the topology of the fibers. In particular, we provide a short review on the existence of the Milnor fibrations in the real and complex cases. This allows us to compare our results with the previous ones. 
\end{abstract}

\begin{section}{Introduction}
In this paper we will work in the category of analytic function-germs. In the book \cite{Milnor} \emph{Singular Points of Complex Hypersurfaces}, John W. Milnor studied singularities of complex hypersurfaces and showed how to attach a locally trivial smooth fibration to each singular point in order to extract topological data from it. Milnor proved the existence of such fibration structures associated to holomorphic functions and also to real analytic maps, but in the latter case for isolated critical points. Both fibrations are currently known as the Milnor fibrations.

\medskip This work is based on five directions: the classical complex case for holomorphic functions (Section~\ref{complex-mf}), 
the real case for analytic map germs with isolated critical points (Section~\ref{real-mf}), the Neuwirth-Stallings (NS) pairs and their relations with polynomial maps with isolated critical point (Section~\ref{NS-pair}) and the topology of Milnor fibrations and their fibers for maps under the Milnor  conditions $(a)$ and $(b)$ as introduced by D. Massey in \cite{Massey} (Section~\ref{tube-fibrations} and Section~\ref{homotopy-type-milnor}).

\medskip We start with a brief review about the existence of Milnor  fibration and its topology in the complex  and real cases (Section~\ref{complex-mf} and Section~\ref{real-mf},  respectively), presenting some old and new recent results (Example~\ref{triv}, Example~\ref{euler-charac} and Theorem~\ref{trivial-fibe}). In Section~\ref{topology-real}  we introduce new results about the topology of the Milnor fibrations (Theorem~\ref{buque}, Proposition~\ref{bouquet-mixed} and Proposition~\ref{not-conclude}). In addition, in Section~\ref{connectness}, we study their fibers with special attention to the real case (Theorem~\ref{Final} and Theorem~\ref{TuboConexo}). 
\end{section}

\begin{section}{Complex Milnor fibrations revisited}\label{complex-mf}

In \cite{Milnor} John Milnor showed that given a non-constant holomorphic function $g: U \subset \mathbb{C}^{n+1} \to \mathbb{C}$ with $g(0) = 0$, there exists a small $\epsilon > 0$ such that

\begin{equation}\label{Fibration}
\phi_{\epsilon} = \displaystyle \frac{g}{\|g\|}: S_{\epsilon}^{2n + 1} \setminus K_{\epsilon} \to S^1
\end{equation}
is a smooth projection of a locally trivial fiber bundle, where $$K_{\epsilon} = g^{-1}(0) \cap S_{\epsilon}^{2n + 1} $$ is called the link of the singularity at origin.

Let $F_{\theta} = \phi_{\epsilon}^{-1} (e^{i\, \theta})$ be the fiber of (\ref{Fibration}) over a point $e^{i\, \theta} \in S^1$. It is a $2n$-dimensional parallelizable real manifold. In addition, using Morse theory, Milnor proved that $F_{\theta}$ has the homotopy type of a finite $CW$-complex of dimension $n$ and the link $K_{\epsilon}$ is $(n-2)$-connected, that is, the homotopy groups $\pi_{j} (K_{\epsilon}) = 0$, for all $j = 0, \ldots, n-2.$

Let $\displaystyle \Sigma_{g} = \{x \in U \subset \mathbb{C}^{n+1} | \nabla g(x) = 0\}$ be the set of critical points of $g$, or the singular locus of $g$. We denote just by $\Sigma$ if no mention to the function is needed.

In the case that $0 \in \Sigma_{g}$ is an isolated point \footnote{In this case we say that $0$ is an \textit{isolated singularity}, or an isolated critical point, of $g$.}, Milnor gave more details about the topology of the fiber and the link. More precisely, in such a case the fiber $F_{\theta}$ has the homotopy type of a wedge of $n$-dimensional spheres $S^n \vee \cdots \vee S^n,$ also known as \emph{Milnor's bouquet of spheres,} with $\mu_{g}$-copies of $S^n$ attached to a single common point. The number $\mu_{g}$ is called \emph{the Milnor number of $g$} at origin. This number is also given by the topological degree of the mapping $$\frac{\bigtriangledown g}{\| \bigtriangledown g \|}: S_{\epsilon}^{2n + 1} \to S^{2n + 1}$$ or, in an algebraic way, it is given by the complex dimension $$\mathrm{dim}_{\mathbb{C}}\, \frac{\mathcal{O}_{n + 1, 0}}{\langle J_g \rangle}$$ where $\mathcal{O}_{n + 1,0}$ is the local ring of germs of holomorphic functions at origin, and $\langle J_g \rangle$ denotes \emph{the Jacobian ideal} which is the ideal spanned by all partial derivatives of first order of $g.$

As a byproduct of the discussions above one easily gets that the Euler-Poincar\'e characteristic of the Milnor fiber is given by

\begin{equation}\label{euler}
\chi(F_{\theta})=1+(-1)^n\mu_{g}
\end{equation}

\begin{exam}\cite{Milnor}
Consider the polynomial
$$g(z_1, z_2) = z_1^3 - z_2^2$$
in two variables, with an isolated critical point at the origin. Then the link
$$K_{\epsilon} = \{(|z_1| e^{2\, i\, \pi t}, |z_2| e^{3\, i\, \pi t}) \, | \, t\in \mathbb{R}, \, |z_1|^2 + |z_2|^2 = \epsilon^2  \}$$
is a trefoil knot in the torus $S_{|z_1|}^1 \times S_{|z_2|}^1$. In this case the fiber $F_{\theta}$ is a Seifert surface with boundary $K_{\epsilon}$ and has the homotopy type of a wedge of $1$-dimensional spheres with $\mu_g = 2.$ Consider Figure \ref{Seifert} as an illustration.
\end{exam}

\begin{figure}[h!]
\centering
\includegraphics[scale=0.35]{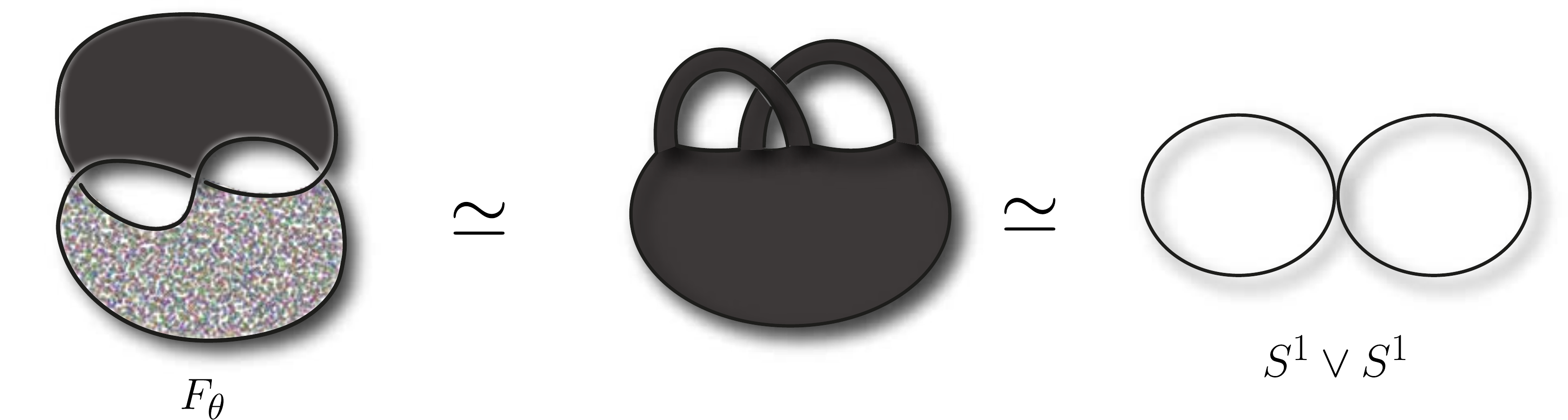}
\caption{The homotopy type of the Milnor fiber.}
\label{Seifert}
\end{figure}

Milnor also proved that for all $\epsilon > 0$ small enough, the manifold
$$(g^{-1} (0) \setminus \{0\}) \cap B_{\epsilon}^{2n +2 }$$
intersects transversally $S_{\epsilon}^{2n + 1}$ and thus $K_{\epsilon}$ is a $(2n-1)$-dimensional smooth manifold. Furthermore, each fiber $F_{\theta}$ can be considered as the interior of a smooth compact manifold with boundary $\overline{F_{\theta}} = F_{\theta} \cup K_{\epsilon}$. Thus in a neighborhood of the link $K_{\epsilon}$ all fibers fit around their common boundary $K_{\epsilon}$ like an open book structure, as illustrated in the Figure \ref{Openbook}.

\begin{figure}[h!]
\centering
\includegraphics[scale=0.24]{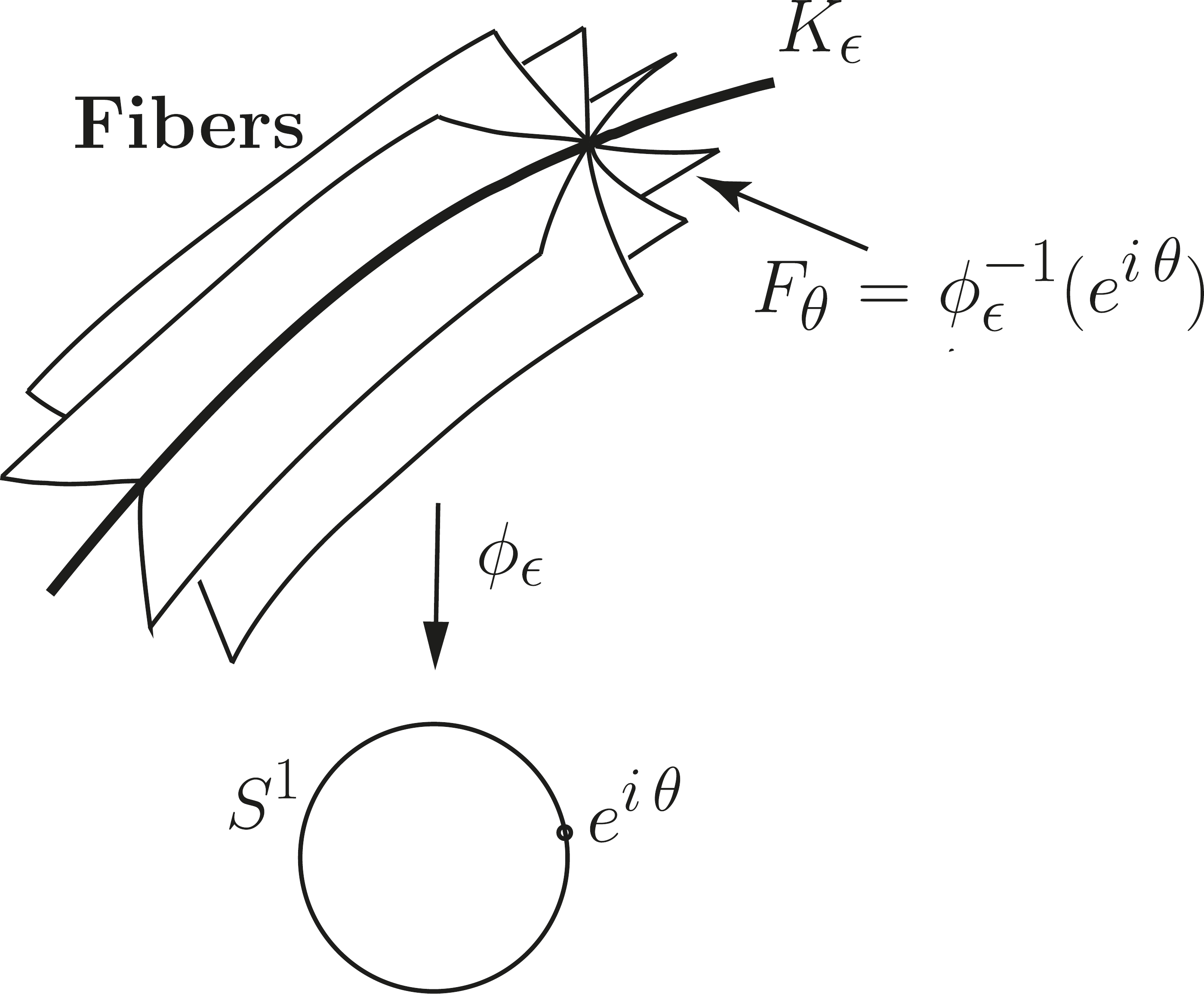}
\caption{Open book structure.}
\label{Openbook}
\end{figure}

Note that in the case of an isolated singular point the Milnor fiber $F_{\theta}$  is $(n-1)$-connected since it has the homotopy type of a wedge of $n$-dimensional spheres $S^n \vee \cdots \vee S^n.$ Later in \cite{KM} Kato and Matsumoto proved an extension of this result as follows below, that in particular shows the close relationship between the singular locus and the topology of the fibers surrounding it.

\begin{theorem} \cite[Theorem 1]{KM} Let $g: U \subset \mathbb{C}^{n+1} \to \mathbb{C}$ be a non-constant holomorphic function defined on an open subset $U$ in the complex $(n+1)$-space $\mathbb{C}^{n+1}$. Suppose that the complex dimension of the critical set $\Sigma_{g}$ is equal to $s$. Then the fiber of Milnor fibration is at least $(n - s - 1)$-connected.
\end{theorem}

\end{section}

\begin{section}{Real Milnor fibrations revisited}\label{real-mf}

Milnor showed the existence of fiber bundle structures for real maps with isolated critical points. More precisely, Milnor considered a real analytic mapping $f: U \subset \mathbb{R}^n \to \mathbb{R}^p$ with $f(0) = 0$ and $n \geq p \geq 2$ such that in some open neighborhood $U$ of the origin $0 \in \mathbb{R}^n$, we have $\Sigma_{f} \cap U = \{ 0 \}$, where $$\Sigma_{f} = \{ x \in U | \mbox{ rank } J f(x) \mbox{ fails to be maximal} \},$$ and $J f(x)$ denotes the Jacobian matrix of $f$ at $x$. This means that $0$ is an isolated singular point of $f$.

We denote by $V=f^{-1}(0)$ and $K_{\epsilon}=V\cap S_{\epsilon}^{n-1}$, where $S_{\epsilon}^{n-1}$ is the sphere of radius $\epsilon$ centered at the origin in $\mathbb{R}^n$. The set $K_{\epsilon}$  is called the \emph{link} of the singularity at the origin. By the transversality theorem one can prove that the diffeomorphism type of $K_{\epsilon}$ does not depend on $\epsilon$ small enough. Hence, we denote $K_{\epsilon}$ just by $K.$

\begin{theorem}\label{TeoMilnorI}\cite[Theorem 11.2]{Milnor} There exists an $\epsilon_0 > 0$ small enough such that for all $\epsilon$ with $0< \epsilon \leq \epsilon_0$, there exists an $\eta_0$ with $0 < \eta_0 \ll \epsilon$ such that the restriction map

\begin{equation}\label{FibraTubo}
f_| : f^{-1} (S^{p-1}_{\eta}) \cap B^{n}_{\epsilon} \to S^{p-1}_{\eta}
\end{equation}
is the projection of a smooth locally trivial fiber bundle for all $\eta$ with
$0 < \eta \leq \eta_0$, where $B^n_{\epsilon}$ denotes the $n$-dimensional closed ball centered at the origin of radius $\epsilon$ in $\mathbb{R}^n$, and $S^{p-1}_\eta$ denotes the sphere of radius $\eta$ centered at the origin in $\mathbb{R}^p$. Furthermore, the fiber $F$ is a smooth compact  $(n-p)-$dimensional manifold bounded  by $K$.
\end{theorem}

The fibration (\ref{FibraTubo}) is known as the {\it real Milnor fibration in the tube} and its fiber $F$ as the {\it real Milnor fiber}. Rescalling the radii $\eta$ by a radial diffeomorphism we can consider the fibration over $S^{p-1},$ the unitary $(p-1)$-dimensional sphere on the target space.

\begin{rema} The real Milnor fibrations have  been extended in several directions for isolated and not isolated singularities. Some interesting results toward that can be found, for example, in \cite{TCR},  \cite{TR}, \cite{Massey} and \cite{RR}.
\end{rema}

An important result about the connectedness of the real Milnor fiber is given in the next proposition, which follows from \cite[Lemma 11.4]{Milnor}.  

\begin{prop} If the link $K = f^{-1}(0) \cap S_{\epsilon}^{n - 1}$ is not empty, then the fiber $F$ of the fibration (\ref{FibraTubo}) is $(p - 2)$-connected.

\end{prop}

In \cite{RDD}, using tools from singularity theory, Morse theory and differential geometry, autors proved an extension of Khimshiashvili's formula (see \cite{Kh}) for the Euler-Poincar\'e characteristic of the real Milnor fiber $F$ of an analytic map germ $f: (\mathbb{R}^n, 0) \to (\mathbb{R}^p, 0)$, $f(x) = (f_1(x), \ldots, f_p(x))$, $n \geq p \geq 2$, with isolated singularity at origin, as below.

\begin{theorem}\label{ThmRDD}\cite{RDD} Let $f$ be an analytic map germ with isolated critical point at origin $0$ and Milnor fiber $F,$ that is the fiber of the associated Milnor fibration (\ref{FibraTubo}).
\begin{itemize}
\item[(a)] If $n$ is even, then  the  Euler-Poincar\'e characteristic of the fiber $F$ is given by $\chi (F) = 1- \mbox{ deg}_0 (\nabla f_1)$, where $\mbox{ deg}_0 (\nabla f_1)$ is the topological degree of the mapping
    $$\frac{\nabla f_1}{\|\nabla f_1\|}: S_{\epsilon}^{n-1} \to S_1^{n-1}.$$
    Moreover, $\mbox{ deg}_0 (\nabla f_1) = \mbox{ deg}_0 (\nabla f_2) = \cdots = \mbox{ deg}_0 (\nabla f_p)$.

\item[(b)] If $n$ is odd, then $\chi (F) = 1$. Moreover $\mbox{ deg}_0 (\nabla f_i) = 0$ for all $i = 1, \ldots, p$.
\end{itemize}
\end{theorem}

\begin{rema}
It is an easy exercise to prove that if we consider a germ of holomorphic function $g:(\mathbb{C}^{n+1},0)\to (\mathbb{C},0)$ with isolated critical point at origin as a real analytic map germ $g:(\mathbb{R}^{2n+2},0)\to (\mathbb{R}^{2},0),$ then the formula $(a)$ above recover the previous Milnor's formula \rm(\ref{euler}\rm). More details can be found in \cite{Andr}.
\end{rema}

\begin{subsection}{Characterization of triviality due to Church and Lamotke}

Due to the Inverse Function Theorem it is known that in a neighbourhood of a non-critical point the fibers of the map behaves like smooth linear planes. However, in the presence of critical points this behaviour is not true in general.

For complex analytic functions $g: U \subset \mathbb{C}^{n+1} \to \mathbb{C},$ $g(0)=0,$ Milnor proved that $g$ has an isolated critical point at origin if and only if the Milnor number $\mu_{g}>0.$ It says that the regular fibers surrounding the singular point must have the homotopy type of a wedge of spheres, and therefore such a singularity cannot be removed even after a topological (continuous) change of coordinates.

However, for real analytic map germ with isolated singularity the picture is not the same in general as we show in the following example.

\begin{exam} \label{triv} Consider $f(x,y,z,w)=\left(x,y(x^2+y^2+z^2+w^2)\right).$ It is easy to see that $\Sigma=\{0\}$ and hence we have a Milnor fibration over $S^{1}_{\eta}$ 
\[f_| : f^{-1} (S^{1}_{\eta}) \cap B^{4}_{\epsilon} \to S^{1}_{\eta}\] The Milnor fiber \begin{align*}
    f^{-1}(\eta,0)&=\{(x,y,z,w)| x=\eta, y(x^2+y^2+z^2+w^2)=0\}\\
    &=\{(\eta,0,z,w)|~z,w\in\mathbb{R}\}
\end{align*} is diffeomorphic to a $2$-plane, and if one consider its intersection with an closed ball centered at origin we get a smooth 2-disk.
\end{exam}

In order to treat this situation Milnor in \cite[p. 100]{Milnor} posed the following question.

\begin{Que}\label{MiQuest}
For which dimensions $n \geq p \geq 2$ do non-trivial examples exist?
\end{Que}
He did not make it clear what does \aspas{non-trivial} should mean. But, for example, he said that certainly the projection $f(x_1, \ldots, x_n) = (x_1, \ldots, x_p)$ is trivial.

However, for a map $f: U \subset \mathbb{R}^n \to \mathbb{R}^p,$ $n\geq p\geq 2,$ with isolated singular point at origin, Milnor proposed the following definition: \aspas{An example will be called trivial if and only if the real Milnor fiber $F$ is diffeomorphic to the closed ball $B^{n-p}$}.
\vspace{0.2cm}

According to his proposal of \aspas{definition} our Example \ref{triv} above is trivial. In a more general way, using Theorem \ref{ThmRDD} and the well known classification of orientable compact connected surfaces with boundary we have the following statement.

\begin{exam}\label{euler-charac} Let $f:U \subset \mathbb{R}^4\to \mathbb{R}^2,$ $0\in U,$ be an analytic map with an isolated critical point at origin and $F$ be the fiber of its Milnor fibration. Then, the map $f$ is trivial if and only if $\chi(F)=1.$ Indeed, the \aspas{if} implication is trivial. The \aspas{only if} implication is a nice exercise, but the reader can see \cite{RDD}, for instance. 
\end{exam}

\vspace{0.2cm}

In \cite{CL} Church and Lamotke answered Milnor's question (Question \ref{MiQuest}) in the following way.

\begin{theorem}\label{Lamotke}\cite[p. 149]{CL}
\begin{itemize}
\item[(a)] For $0 \leq n - p \leq 2$, non-trivial examples occur precisely for the dimensions $(n, p) \in \{ (2,2), (4,3), (4, 2)\}$.

\item[(b)] For $n - p \geq 4$, non-trivial examples occur for all $(n, p)$.

\item[(c)] For $n - p = 3$, all examples are trivial except $(5, 2)$, $(8, 5)$ and possibly $(6, 3)$.
\end{itemize}
\end{theorem}

Furthermore, they also presented an alternative characterization of triviality in terms of the branch set as we describe below.

\vspace{0.2cm}

A map $f: X \to Y$ is \emph{locally topologically equivalent} at $x\in X$ to a map $g: X^{\prime} \to Y^{\prime}$ at $x^{\prime} \in X^{\prime}$ if there are open neighborhoods $U$ of $x$, $U^{\prime}$ of $x^{\prime}$, $V$ of $f(x)$, $V^{\prime}$ of $g(x^{\prime})$ and homeomorphisms $\alpha: U \to U^{\prime}$ and $\beta: V \to V^{\prime}$ such that
$$\beta \circ f|_{U} = g \circ \alpha.$$

\begin{defi} The branch set of $f$, denoted by $B_f$, is the set of points at which $f$ fails to be  locally topologically equivalent to the (canonical) projection $\mathbb{R}^n \to \mathbb{R}^p$.
\end{defi}

By the Inverse Function Theorem we have that $B_{f} \subseteq \Sigma .$ However, as Example \ref{triv} shows the set $B_{f}$ may be empty whereas $\Sigma$ may not be.

\begin{prop}\label{CeL}\cite[p. 151]{CL} For $n - p \neq 4, \, 5$, a real polynomial map $f: (\mathbb{R}^n, 0) \to (\mathbb{R}^p, 0)$, with $0$ an isolated critical point, is trivial if and only if $0 \not \in B_f$.
\end{prop}

\begin{rema}
In the proof of Proposition \ref{CeL},  Church and Lamotke used the Poincar\'e conjecture in dimensions $ n - p -1 \leq 2$. Since the Poincar\'e conjecture is true in the dimension 3 (proved by Perelman), then we can use the same argument as the one used by Church and Lamotke to extend the result for the dimension $ n - p - 1 = 3$, that is, $n - p = 4$.
\end{rema}

It follows from Theorem \ref{Lamotke} that in the case $n-p = 3$ the map $f$ can be nontrivial only if $(n, p) \in \{(5,2),\, (6,3), \, (8,5)\}$.

In the theorem below we give conditions under which the singular map germ is trivial for the pair $(5,2)$.

\begin{theorem}\label{trivial-fibe}
 Let $f: (\mathbb{R}^5, 0) \to (\mathbb{R}^2, 0)$ be an  analytic map germ with an isolated singularity at origin. Then, $f$ is trivial if and only if the link $K$ is an unknotted $2$-sphere in the $4$-dimensional sphere $S^4$.
\end{theorem}

\proof If $f$ is trivial then the link $K$, which is the boundary of the Milnor fiber $F$, must be an unknotted sphere.

Conversely, suppose that $K$ is isotopic to the $2$-sphere $S^2$. According to \cite[Lemma 11.3]{Milnor} the total space $E$ of the Milnor fibration is diffeomorphic to the complement of an open tubular neighborhood of $K$ in $S^4$. Therefore, $E$ has the same homotopy type of $S^4 \setminus (S^2 \times B^2)$. 

We have $S^4 \setminus (S^2 \times B^2) = B^3 \times S^1$. In fact
$$\begin{array}{lll}
   S^4  = \partial B^5 &  = & \partial (B^3 \times B^2)\\
   & = & (\partial B^3 \times B^2) \cup (B^3 \times \partial B^2)\\
   &=& (S^2 \times B^2) \cup (B^3 \times S^1).\\
\end{array}$$
Thus $S^4 \setminus (S^2 \times B^2) =  (S^2 \times B^2) \cup (B^3 \times S^1) \setminus (S^2 \times B^2) = B^3 \times S^1.$

Consider the following piece of the long exact sequence in homotopy of the Milnor fibration $F \hookrightarrow E \rightarrow S^1$:
\begin{equation}\label{Sequencia}
  \pi_2 (S^1) \to \pi_1 (F) \to \pi_1 (E) \to \pi_1 (S^1) \to \pi_0 (F).
\end{equation}

We have $\pi_0 (F) = 0$ and $\pi_1 (E) = \pi_1 (B^3 \times S^1) = \pi_1(S^1) = \mathbb{Z}$. Thus we can write (\ref{Sequencia}) as follows:

$$0 \to \pi_1 (F) \to \mathbb{Z} \to \mathbb{Z} \to 0.$$
Then $\pi_1 (F) = 0$ and therefore $F$ is simply connected.

Let $N = F \cup_{K} F$ be the double of $F$  that is a closed orientable $3$-manifold. We have, by the Van Kampen Theorem,
$$\pi_1 (N) \cong \pi_1(F)\ast\pi_1(F).$$
Therefore $\pi_1 (N) = 0$. Thus, by the Poincar\'e conjecture (proved by Perelman) we conclude that $N$ is diffeomorphic to a $3$-sphere and consequently $F$ is diffeomorphic to a $3$-disk, that is, $f$ is trivial. \cqd
\end{subsection}
\end{section}

\begin{section}{The Milnor fibration from $NS$-pairs}\label{NS-pair}

As a consequence of  Theorem \ref{TeoMilnorI}, Milnor used an appropriate vector fields to prove the following statement.

\begin{theorem}\label{TeoMilnorII} Let $f: (\mathbb{R}^n, 0) \to (\mathbb{R}^p, 0)$, $n \geq p \geq 2$, be an analytic map germ with an isolated singularity at the origin. Then, there exists an $\epsilon_0 > 0$ such that for all $0 < \epsilon \leq \epsilon_0$, the complement of the link $K = f^{-1} (0) \cap S_{\epsilon}^{n-1}$ in  $S_{\epsilon}^{n-1}$ is the total space of a smooth fiber bundle over the sphere $S^{p-1}$ with each fiber $F_{f}$ being the interior of a compact $(n - p)$-dimensional manifold bounded by a copy of $K$, where $S_{\epsilon}^{n-1}$ denotes the sphere in $\mathbb{R}^n$ with radius $\epsilon$ centred at the origin.
\end{theorem}


Motivated by Theorem \ref{TeoMilnorII}, Looijenga in \cite{Looijenga} formulated the following notion.

\begin{defi} \cite[p. 418]{Looijenga} Let $K = K^{n-p-1}$ be an oriented submanifold of dimension $n - p - 1$ of the oriented sphere $S^{n-1}$ with trivial normal bundle, or let $K = \emptyset$. Suppose that for some trivialization $c: N(K) \to K \times D^p$ of a tubular neighborhood $N(K)$ of $K$ the fiber bundle defined by the composition
$$N(K) \setminus K \stackrel{c}{\to} K \times (D^p \setminus \{ 0 \}) \stackrel{\pi}{\to} S^{p-1}, 
$$ with the last projection given by $\pi (x, y) = y/ \|y \|$, extends to a smooth fiber bundle $S^{n-1} \setminus K \to S^{p-1}$. Then the pair $\left(S^{n - 1}, K^{n - p - 1}\right)$ is called a \emph{Neuwirth-Stallings pair}, or an $NS$-pair for short.
\end{defi}

It follows from the proof of Theorem \ref{TeoMilnorII} and a \aspas{little more effort} that $\left(S_{\epsilon}^{n-1}, f^{-1} (0) \cap S_{\epsilon}^{n-1}\right)$ is an $NS$-pair. In this case, it is called  the $NS$-pair associated to the singularity.

Looijenga in \cite[pp. 419-420]{Looijenga} showed how to use the connected sum of $NS$-pairs to construct new $NS$-pairs but with an additional property. In fact, he proved that given an $NS$-pair $\left(S^{n - 1}, K^{n - p - 1}\right)$ with fiber $F$, after performing a special connected sum of itself, that he denoted by $$\left(S^{n - 1}, K^{n - p - 1}\right) \sharp \left((-1)^{n-1} S^{n - 1}, (-1)^{n-p} K^{n - p - 1}\right),$$ there exists a polynomial map germ $f: (\mathbb{R}^{n}, 0) \to (\mathbb{R}^{p}, 0)$ with an isolated singularity at the origin such that the associated $NS$-pair of $f$ (or, its Milnor fibration on the sphere) is isomorphic to that connected sum, with the fiber being diffeomorphic to the interior of $\overline{F} \natural (-1)^{n-p} \overline{F}$, where $\natural $ means the connected sum along the boundary $\partial F$.

\subsection{Updating Milnor's question in the (6,3) case}

Note that, according to Theorem \ref{Lamotke} a polynomial map germ $f: (\mathbb{R}^n, 0) \to (\mathbb{R}^p, 0)$, $n \geq p \geq 2$  with an isolated singularity at the origin, with $n - p = 3,$ may be non-trivial only if $(n, p) \in \{ (6, 3), (8, 5), (5, 2)\}.$ Therefore, after Church-Lamotke's Theorem the Milnor question \ref{MiQuest} remained open only in the pair of dimension $(6,3).$

In \cite[Section 3]{Saeki}, the authors used the topology of configurations spaces to prove a characterizations of Neuwirth-Stallings pair $(S^5, K)$ with \linebreak $2$-dimensional $K$. From this characterization and Looijenga's construction described above they showed the existence of polynomials map germs $(\mathbb{R}^6, 0) \to (\mathbb{R}^3, 0),$ with an isolated singularity at origin, such that the associated Milnor fibers are not diffeomorphic to a disk. Thus putting an end to Milnor's question (Question \ref{MiQuest}). Further details can be seen in \cite[Proposition 3.8]{Saeki}.

\end{section}

\begin{section}{Tube fibrations for non-isolated singularities}\label{tube-fibrations}

In this section we will follow the definitions and results given by D. Massey in the paper \cite{Massey}.

Let $f = (f_1, \ldots, f_p): (\mathbb{R}^n, 0) \to (\mathbb{R}^p, 0)$ be a non-constant analytic map, $1\leq p \leq n-1$, $V:= f^{-1}(0)$ and $\Sigma_{f}$ be the set of critical points of $f.$

Let $r: \mathbb{R}^n \to \mathbb{R}$ be the square of the distance function to the origin, $r(x)=\|x\|^2$, and let $\Sigma_{(f, r)}$ be the set of critical points of the pair $(f,r)$, that is, the set of points where the gradients $\nabla r,  \nabla f_1, \ldots, \nabla f_p$ are linearly dependent. Note that $\Sigma_{f} \subseteq \Sigma_{(f,r)}.$

\begin{defi}\cite{Massey} Given $f$ and $r$ as above.

\begin{itemize}

\item[(i)] The map $f$ satisfies the \emph{Milnor condition $(a)$ at origin} if $\Sigma_{f} \subset V$ in a neighborhood of the origin.

\item[(ii)] The map $f$ satisfies the \emph{Milnor condition $(b)$ at origin} if $0$ is isolated in $V\cap \overline{\Sigma_{(f, r)}\setminus V}$ in a neighborhood of the origin, where the \aspas{bar} notation means the topological closure of the space $\Sigma_{(f, r)}\setminus V$ in $\Sigma_{(f,r)}.$

\item[(iii)] If $f$ satifies the Milnor's condition $(a)$ \rm(respectively $(b)$\rm), then we say that $\epsilon_{0} >0$ is a \emph{Milnor $(a)$ radius for $f$ at origin} (respectively, Milnor $(b)$ radius for $f$ at origin) provided that for all $0<\epsilon \leq \epsilon_{0}$ one gets
$B_{\epsilon}^n \cap (\overline{\Sigma_{f} \setminus V}) = \emptyset$ \rm(respectively, $B_{\epsilon}^n \cap V \cap  (\overline{\Sigma_{(f, r)} \setminus V}) \subseteq \{0\}$\rm).
\end{itemize}
\end{defi}

Under condition $(iii)$ above, we simply say that $\epsilon_{0}$ is a \emph{Milnor radius for $f$ at origin}, if $\epsilon_{0}$ is both a Milnor $(a)$ and Milnor $(b)$ radius for $f$ at origin.

\begin{rema} It can be concluded that, under Milnor's conditions $(a)$ and $(b)$, for all regular values close to the origin of $\mathbb{R}^p$, the respective fibers inside the closed $\epsilon$-ball are smooth and transverse to the sphere $S_{\epsilon}^{n-1}$.
\end{rema}

\begin{rema} It is easy to see that any analytic map germ $f = (f_1, \ldots, f_p): (\mathbb{R}^n, 0) \to (\mathbb{R}^p, 0)$ with isolated critical point at origin has a Milnor's radius. Indeed, condition (a) is satisfied because, by assumption, the critical point is isolated. Condition (b) is satisfied because $V\setminus\{0\}$ is transversal to small spheres, and so are small levels of $f$.

\end{rema}

In a more general way, the Milnor conditions $(a)$ and $(b)$ are enough to ensure the existence of the Milnor tube fibrations as follows.

\begin{theorem}\label{mtu}\cite[p. 284, Theorem 4.3]{Massey} Let $f = (f_1, \ldots, f_p): (\mathbb{R}^n, 0) \to (\mathbb{R}^p, 0)$ be an analytic map germ, and $\epsilon_0$ be a Milnor radius for $f$ at the origin. Then, for each $0< \epsilon \leq \epsilon_0$, there exists $\eta$, $0< \eta \ll \epsilon$, such that

\begin{equation}\label{MasseyI}
f_{|}: B_{\epsilon}^n \cap f^{-1}(B_{\eta}^p \setminus \{0\}) \to B_{\eta}^p \setminus \{0\},
\end{equation}
is the projection of a smooth locally trivial fibration.
\end{theorem}

\begin{cor} {\rm(Fibration in the Milnor tube)} Let $f = (f_1, \ldots, f_p): (\mathbb{R}^n, 0) \to (\mathbb{R}^p, 0)$ be an analytic map germ, and $\epsilon_0$ be a Milnor radius for $f$ at the origin. Then, for each $0< \epsilon \leq \epsilon_0$, there exists $\eta$, $0< \eta \ll \epsilon$, such that

\begin{equation}\label{MasseyII}
f_{|}: B_{\epsilon}^n \cap f^{-1}(S_{\eta}^{p-1}) \to S_{\eta}^{p-1},
\end{equation}
is the projection of a smooth locally trivial fibration.
\end{cor}

Given $f = (f_1, \ldots, f_p): (\mathbb{R}^n, 0) \to (\mathbb{R}^p, 0),$ $n>p\geq 2,$ a non-constant analytic map germ and $V= f^{-1}(0),$ H. Whitney proved that one can find a stratification of $V,$ let say $\mathcal{S},$ with finitely many subanalytics strata such that each stratum $S_{i}\in \mathcal{S}$ that contains the origin $\{0\}$ on its closure is transversal to all sphere $S_{\epsilon}^{n-1},$ for all $\epsilon>0$ small enough. Such a stratification induces a stratification of the link $K_{\epsilon}=V\cap S_{\epsilon}^{n-1}.$ Moreover, one can use such transversality to prove that the topological type of $K_{\epsilon}$ does not depend on $\epsilon$ small enough. For this reason from now on we will denote the link only by $K.$
\end{section}

\begin{section}{Homotopy type of the Milnor fibers revisited}\label{homotopy-type-milnor}

Let $f = (f_1, \ldots, f_p): (\mathbb{R}^n, 0) \to (\mathbb{R}^p, 0)$ be an analytic map germ with isolated singularity and set $\phi = (f_1, \ldots, f_{p - 1})$ be the composition of $f$ with the canonical projection $\Pi: \mathbb{R}^p \to \mathbb{R}^{p -1}, \Pi(y_{1},\ldots,y_{p})=(y_{1},\ldots,y_{p-1}).$

In \cite[p. 100]{Milnor} Milnor conjectured that the fiber of $\phi$ is homeomorphic to the product of the unit interval with the fiber associated to $f$. This conjecture was proved true in more generally in \cite{RD} as described below.

\begin{lema}\cite[Lemma 4.1]{RD} If $f$ satisfies Milnor's conditions $(a)$ and $(b)$, then $\phi=\Pi\circ f$ also satisfies Milnor's conditions $(a)$ and $(b)$.
\end{lema}

Let $F_{f}$ be the Milnor fiber of the map $f$ and $F_{\phi}$ be the Milnor fiber of its projection $\phi,$ regarding the Milnor tube fibration \ref{MasseyII}. In \cite{RD} the authors proved the following.

\begin{theorem}\cite[Theorem 6.3]{RD} Let $f = (f_1, \ldots, f_p): (\mathbb{R}^n, 0) \to (\mathbb{R}^p, 0)$, $p \geq 2$, be  an analytic mapping that satisfies Milnor's conditions $(a)$ and $(b)$ and let $\phi=\Pi \circ f: (\mathbb{R}^n, 0) \to (\mathbb{R}^{p-1}, 0)$. Then, the fiber $F_{\phi}$ is homeomorphic to $F_f \times [-1, 1]$.
\end{theorem}

\begin{cor}\cite[Milnor's conjecture]{RD} Let $f = (f_1, \ldots, f_p): (\mathbb{R}^n, 0) \to (\mathbb{R}^p, 0)$, $p \geq 2$, be  an analytic mapping with isolated singular point at origin, and $\phi=\Pi \circ f: (\mathbb{R}^n,0) \to (\mathbb{R}^{p-1}, 0)$. Then, the Milnor fiber $F_{\phi}$ is homeomorphic to $F_f \times [-1, 1].$
\end{cor}

\end{section}

\begin{section}{Topology of the real Milnor fibers}\label{topology-real}

In this section we prove some improvements and new results about the topology of the real Milnor fibers. 

In the proofs of this section we will consider the homology and cohomology groups with integer coefficients. For our purpose we present and prove the following statement.  

\begin{lema}\label{torsionfree}
Let $M$ be a $m$-dimensional manifold ($m \geq 3$), compact, orientable, with boundary and simply connected. Then the homology groups  $H_{m-1}(M)$ and $H_{m-2}(M)$ are torsion free.
\end{lema}

\proof We have that $(H_i (M)/T_i) \oplus T_{i-1} \cong H^i (M)$  where $T_{i}$ is the torsion subgroup of $H^{i} (M)$. Then, $H_{m-1} (M)$ is torsion free, since $H^m (M) = 0$.

On the other hand, by duality
$$H^{m-1} (M) = H_{m - (m-1)} (M, \partial M) = H_1 (M, \partial M).$$

Consider the following piece of the long exact sequence in homology of the pair $(M, \partial M)$

$$0 \to H_1 (M, \partial M) \to H_0 (\partial M) \to H_0 (M) \to H_0 (M, \partial M) = 0,$$
then 
$$0 \to H_1 (M, \partial M) \to \displaystyle\oplus_{\alpha} \mathbb{Z} \to \mathbb{Z} \to 0,$$
where $\alpha$ is the number of the connected components of $\partial M$. Consequently
$$H_1(M, \partial M) \oplus \mathbb{Z} = \oplus_{\alpha} \mathbb{Z}$$
and therefore $T_{m-2}$ is trivial, that is, $H_{m-2} (M)$ is torsion free.\cqd

We are ready to state and prove the following theorem. 

\begin{theorem}\label{buque}Let $f: (\mathbb{R}^{2p + 1}, 0) \to (\mathbb{R}^p, 0)$ be a non-trivial real analytic map with an isolated singularity at the origin. If $p\geq 3$ and $f$ is not trivial, then the Milnor fiber $F$ has the homotopy type of a bouquet of spheres of the form
$$\bigvee_{i=1}^{\beta_{p-1}} (S^{p-1} \vee S^p ),$$
where $\beta_{p-1}$ is $(p - 1)$-th Betti number of the fiber $F$.
\end{theorem}

\proof  
The fiber $F$ is a $(p+1)$-dimensional manifold, compact, orientable, with boundary and $(p-2)$-connected then, by Lemma \ref{torsionfree}, $H_{p} (F)$ and $H_{p-1} (F)$ are torsion free.

By the Hurewicz theorem, for $1\leq j \leq p-2$,
$$\pi_j (F) = H_j (F) = 0$$
and the Hurewicz homomorphisms
$$\rho_{p-1}: \pi_{p-1} (F) \to H_{p-1} (F) = \mathbb{Z}^{\beta_{p-1}},$$
$$\rho_{p}: \pi_{p} (F) \to H_{p} (F) = \mathbb{Z}^{\beta_{p}}$$
are surjective.

Note that by Theorem \ref{ThmRDD} we have that $$\beta_{p-1} = \beta_p.$$

On the other hand, $\rho_j ([\psi]) = \psi_{*} ([S^j])$, where $[S^j]$ is a generator of $H_j (S^j) = \mathbb{Z}$. Thus, for each generator
$\gamma_i \in H_{p-1} (F) = \mathbb{Z}^{\beta_{p-1}}$, $1\leq i \leq \beta_{p-1}$,  there is $\psi_i: S^{p-1} \to F$ such that $\rho_{p-1} ([\psi_i]) = \gamma_i$. Similarly, for each generator $\Gamma_i \in H_{p} (F) = \mathbb{Z}^{\beta_{p}}$, $1\leq i \leq \beta_{p-1}$,  there is $\eta_i: S^p \to F$ such that $\rho_{p} ([\eta_i]) = \Gamma_i$. Therefore, we have the continuous map
$$l= \bigvee_{i=1}^{\beta_{p-1}} (\psi_i \vee \eta_i): \bigvee_{i=1}^{\beta_{p-1}} (S^{p-1} \vee S^p ) \to F,$$
obtained by the wedge of the maps $\psi_i$ and $\eta_i$, for each $i = 1, 2, \ldots, \beta_{p-1}$, which is a homotopy equivalence by the Whitehead theorem.
\cqd

\vspace{0.2cm} From \cite{Ok} we recall that a \emph{mixed polynomial} $f$ of $n$-variables $z=(z_1, \ldots, z_n)$ in $\mathbb{C}^n$ is a function expanded in a convergent power series of variables $z=(z_1, \ldots, z_n)$ and $\overline{z}=(\overline{z}_1, \ldots, \overline{z}_n)$, that is,
$$f(z, \overline{z}) = \sum_{\nu, \mu} c_{\nu, \mu} z^{\nu}\overline{z}^{\mu} $$
where $z^{\nu} = (z_1^{{\nu}_1}, \ldots, z_n^{{\nu}_n})$ for $\nu = (\nu_1, \ldots, \nu_n)$ (respectively $\overline{z}^{\mu}=(\overline{z}_1^{{\mu}_1}, \ldots, \overline{z}_n^{{\mu}_n})$ for $\mu = (\mu_1, \ldots, \mu_n)$) as usual. Here $\overline{z}_j$ is the complex conjugate of $z_j$.

\begin{defi}\cite{Ok2} A mixed polynomial $\displaystyle f(z, \overline{z}) = \sum_{\nu, \mu} c_{\nu, \mu} z^{\nu}\overline{z}^{\mu}$ satisfies the  Hamm-L\^e condition at the origin if there exists a positive number $r_0$ such that for any $0< r_1 \leq r_0$,  there exists a positive number $\delta(r_1)$ such that the hypersurface $f^{-1}(\eta) \cap B^{2n}_{r_0}$ is non-empty, non-singular and it intersects transversely with the sphere $S^{2n-1}_{r}$ for any $r, \eta$ with $r_1\leq r \leq r_0$ and $0\neq |\eta| \leq \delta(r_1)$. Here $ B^{2n}_{r}$ is the ball $\{z\in \mathbb{C}^n \, | \,  \|z\|\leq r\}$ of radius $r$ and $S^{2n-1}_{r}$ is the boundary sphere of $ B^{2n}_{r}$. We call such a positive number $r_0$ a stable Milnor radius.
\end{defi}

\begin{prop}\label{bouquet-mixed}
Assume that $f(z, \overline{z})$ is a mixed function of $3$-complex variables which satisfies the  Hamm-L\^e condition at the origin. Consider a stable Milnor radius $r_0$ and assume that there exists a mixed smooth point $w \in f^{-1} (0)$ with $\|w\| < r_0$ and the sphere of radius  $\|w\|$ intersects $f^{-1} (0)$  transversely at $w$. We assume also that $\mathrm{dim}_{\mathbb{R}} f^{-1} (0) = 4$.
If  $\pi_1 (F)$ is trivial then $F$ has the homotopy type of a bouquet 
$$\bigvee_{i=1}^{\beta_2} S^2\vee \bigvee_{i=1}^{\beta_3}S^3,$$
where $\beta_2$ and $\beta_3$ are, 
respectively, the second  Betti-number and the third Betti-number of the fiber $F$.
\end{prop}

\proof By \cite{Ok2} the Milnor fiber $F$ is connected and by hypothesis $\pi_1 (F)$ is trivial.
Thus $F$ is a $4$-dimensional, simply connected, oriented, compact manifold with boundary then, by Lemma \ref{torsionfree}, $H_3 (F)$ and $H_2 (F)$ are torsion free. Therefore, using an argument similar to that used in Theorem \ref{buque} we conclude that the Milnor fiber has the same homotopy type of a bouquet of spheres

$$\bigvee_{i=1}^{\beta_2} S^2\vee \bigvee_{i=1}^{\beta_3}S^3.$$
\cqd

\vspace{0.2cm}

\begin{rema} An analogous result to Theorem \ref{buque} is proved in \cite{Saeki} for $n = 2p$. More precisely, let $f: (\mathbb{R}^{2p}, 0) \to (\mathbb{R}^p, 0)$ be a non-trivial real analytic map with an isolated singularity at the origin then the Milnor fiber is a
$(p - 2)$-connected $p$-dimensional manifold with boundary. Therefore $H_{p-1} (F)$ is torsion free and by an argument similar to that used in Theorem \ref{buque} we can conclude that the Milnor fiber has the same homotopy type of a bouquet of spheres $\displaystyle \bigvee_{i=1}^{\beta_{p-1}} S^{p-1}$.
\end{rema}

\vspace{0.2cm}

The following proposition shows that we can not expect, in general, a bouquet Theorem in the real case like in the complex case for an isolated critical point.

\begin{figure}[h!]
\centering
\includegraphics[scale=0.25]{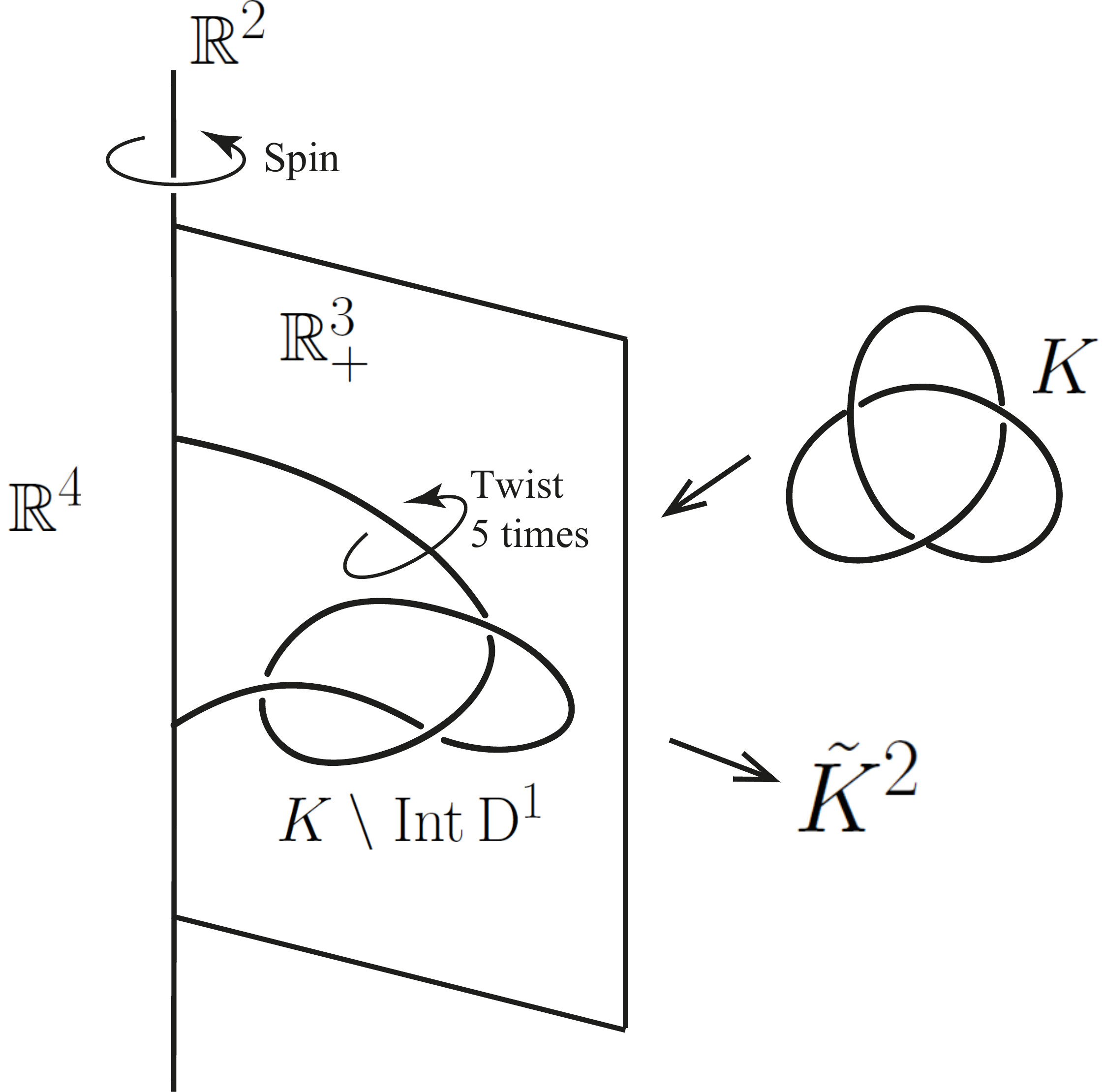}
\caption{Twist-spinning.}
\label{twist-spun}
\end{figure}

\begin{prop}\label{not-conclude} There exists a real polynomial map with an isolated singularity at the origin $f: (\mathbb{R}^5, 0) \to (\mathbb{R}^2, 0)$ non trivial such that the corresponding Milnor fiber does not have the homotopy type of a bouquet of spheres.
\end{prop}

\proof   Using a process due to Zeeman \cite{Zeeman} known as \aspas{twist-spinning} we obtain an $NS$-pair $(S^4, \tilde{K}^2)$, being $\tilde{K}^2$ the 5-twist-spun trefoil.  Figure~\ref{twist-spun} illustrates this construction.

According to Zeeman, the fiber $F$ of this $NS$-pair is such that
$$\pi_1 (\overline{F}) = \left< a, \, b, \, c\,  | \, a^5 = b^3 = c^2 = abc \right>  \mbox{ (the dodecahedral group)}.$$
Then Looijenga's construction leads to non-trivial polynomial map $f: (\mathbb{R}^5, 0) \to (\mathbb{R}^2, 0)$ with an isolated singularity at origin that has as fiber the interior of the connected sum along the boundary $\overline{F} \natural -\overline{F}$. However,
$$\pi_1 (\overline{F} \natural -\overline{F}) = \pi_1 (\overline{F}) \ast \pi_1 (-\overline{F}) \mbox{ (free product)},$$
which does not correspond to the fundamental group of a bouquet of \linebreak spheres. \cqd
\end{section}

\section{Connectness of the Milnor tube and fibers}\label{connectness}

In this section we introduce new results about the topology of the Milnor fibrations and fibers. For that we start with the following statement.

\begin{theorem}\label{Final} Let $f = (f_1, \ldots, f_p): (\mathbb{R}^n, 0) \to (\mathbb{R}^p, 0)$ be an analytic map that satisfies Milnor's conditions $(a)$ and $(b)$ at origin. Then the total space $E = B_{\epsilon}^n \cap f^{-1} (S_{\delta}^{p-1})$ of (\ref{MasseyII}) is diffeomorphic to the complement $S_{\epsilon}^{n-1} \setminus T$, where $T$ is a subanalytic neighborhood of $K:= f^{-1} (0) \cap S_{\epsilon}^{n-1}$ in $S_{\epsilon}^{n-1}$.
\end{theorem}

\proof From \cite[Corollary 1.3 and Proposition 1.6]{Durfee} we consider the subanalytic set 
$$T = \{x \in S_{\epsilon}^{n-1} |~ 0\leq \alpha (x) \leq \gamma\}$$ where $\alpha (x)$ is given by 
$$\alpha (x) = \|f(x)\|^{2}$$
with $\gamma$ smaller than all nonzero critical values of $\alpha$. Furthermore, $K = \alpha^{-1}(0)\cap S_{\epsilon}^{n-1}$.

Consider the vector field $v(x) := \nabla \alpha (x)$ on $B_{\epsilon}^n \setminus V$. The euclidean inner products $\langle v(x), x\rangle$ and
$\langle v(x), \nabla \|f(x)\|^{2}\rangle$ are both positive since $\nabla \|f(x)\|^{2}$ and $2 x = \nabla \|x\|^{2}$ are nonzero throughout $B_{\epsilon}^n \setminus V$ and cannot point in opposite directions by \cite[Corollary 3.4]{Milnor}.

\begin{figure}[!ht]
\centering
\includegraphics[scale=0.2]{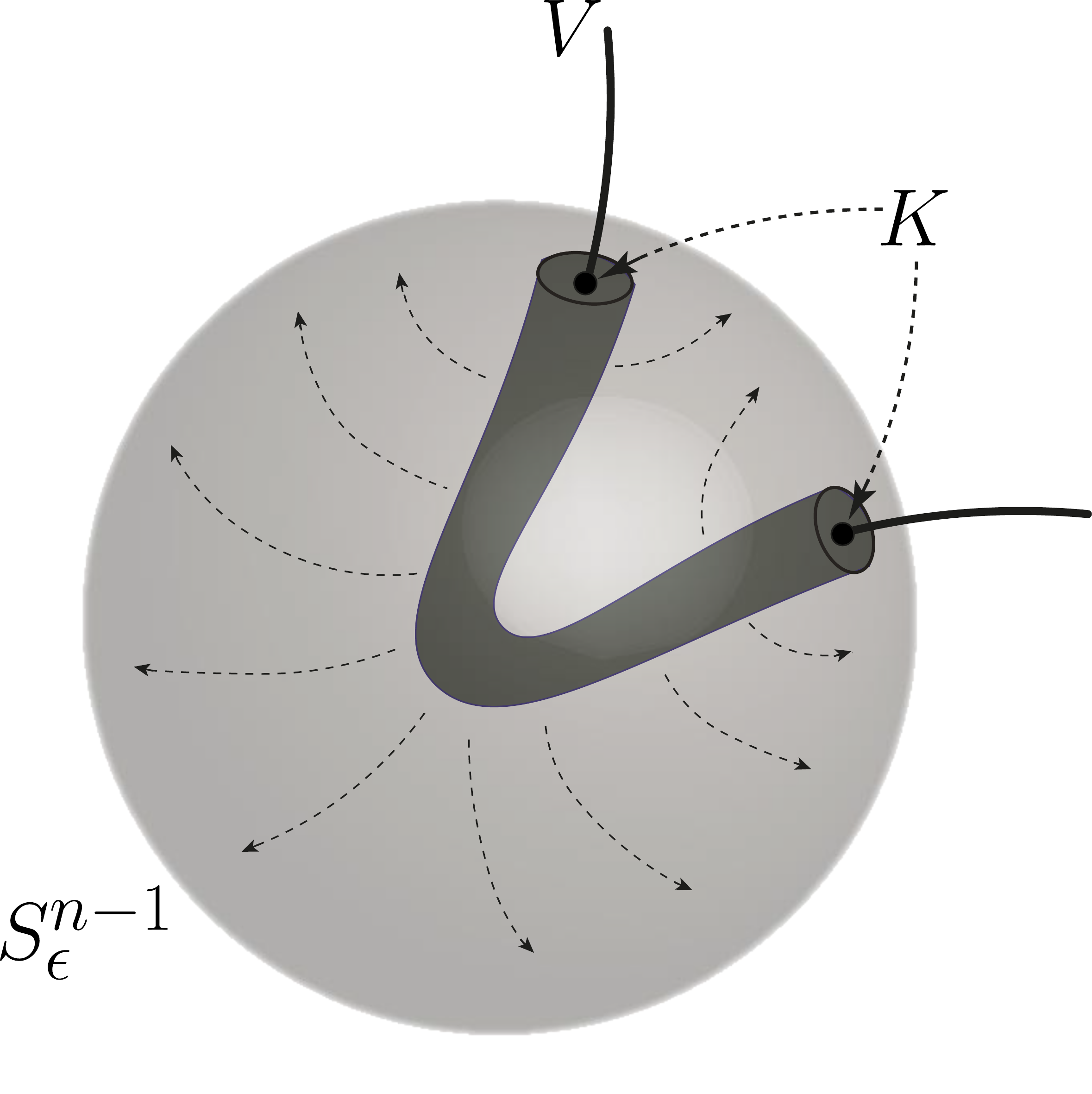}
\caption{Vector field.}
\label{tubo}
\end{figure}

Now pushing out along the trajectories of the vector field $v$ we obtain a diffeomorphism from $B_{\epsilon}^n \cap f^{-1}(S_{\delta}^{p-1}) $ to $S_{\epsilon}^{n-1} \setminus T$, see  Figure~\ref{tubo} as an illustration. More precisely, we start at any point $z$ of $B_{\epsilon}^n \cap f^{-1} (S_{\delta}^{p-1})$ and follow the trajectory  of $v$ through $z$  until we reach a point $z'$ on $S_{\epsilon}^{n-1} \setminus T$. The correspondence $z\mapsto z'$ is a diffeomorphism from $B_{\epsilon}^n \cap f^{-1}(S_{\delta}^{p-1}) $ to $S_{\epsilon}^{n-1} \setminus T$.\cqd

\begin{cor} If $f: (\mathbb{R}^n, 0) \to (\mathbb{R}^{2}, 0)$ is an analytic map that satisfies Milnor's conditions $(a)$ and $(b)$ at origin then
$\pi_i(F) = \pi_i (S_{\epsilon}^{n-1} \setminus T),$
for  all $i\geq 2$.
\end{cor}

\proof Consider the following piece of the long exact sequence in homotopy of the Milnor fibration: 
$$\pi_{i+1}(S^{1}) \to \pi_i (F) \to \pi_i (E) \to \pi_i(S^1).$$
If $i\geq 2$ then $\pi_{i+1}(S^{1}) = \pi_{i}(S^{1})  = 0$ and therefore, $\pi_i (F) = \pi_i (E) = \pi_i (S_{\epsilon}^{n-1} \setminus T)$.
\cqd

The next contribution of this section is the following statement. 

\begin{theorem}\label{TuboConexo}
Let $f = (f_1, \ldots, f_p): (\mathbb{R}^n, 0) \to (\mathbb{R}^p, 0)$, $n > p \geq 2$, be an analytic map that satisfies Milnor's conditions $(a)$ and $(b)$ at origin. Assume further that all cells of $K$ have dimensions $\leq n-p-1.$

\begin{itemize}
\item[(i)] If $p\geq 2,$ then the total space of the Milnor fibration (\ref{MasseyII}) is path connected.

\item[(ii)] If $p \geq 3$ then the Milnor fiber is path connected.

\item[(iii)] If $p=2$ the Milnor fiber is path connected if and only if the Milnor fibration (\ref{MasseyII}) admits a (global) cross-section.

\end{itemize}

\end{theorem}

\proof $(i)$ Let $E$ denote the total space of (\ref{MasseyII}). By Theorem \ref{Final}, we have
$$\widetilde{H}_q (E) = \widetilde{H}_q (S_{\epsilon}^{n-1} \setminus T) \mbox{ (reduced homology)}$$
According to \cite{lojasiewicz1964}, a subanalytic set can be triangulated and therefore is locally contractible. Then, by Alexander duality,
$$\widetilde{H}_q (S_{\epsilon}^{n-1} \setminus T) = \widetilde{H}^{n -2 - q} (T).$$

Note that if $p - q >1$ then $\widetilde{H}^{n - 2 - q} (K) = 0$. Hence, since $n > p \geq 2$ and $K \hookrightarrow T$ is a homotopy equivalence, we have
$$\widetilde{H}_0 (E) = \widetilde{H}_0 (S_{\epsilon}^{n-1} \setminus T) = \widetilde{H}^{n-2}(T) = \widetilde{H}^{n-2}(K) = 0,$$
because $n-2>n-p-1$ and we are assuming that all cells of $K$ have dimension $ \leq n-p-1.$ Thus $H_0 (E) = \mathbb{Z}$ and therefore $E$ is path connected.

 \vspace{0.3cm}

\noindent $(ii)$ Consider the following piece of the long exact sequence in homotopy of the Milnor fibration: 
$$\pi_1(S^{p-1}) \to \pi_0 (F) \to \pi_0 (E).$$
Note that $\pi_1(S^{p-1}) = 0$, since $p-1\geq 2$, and $\pi_0 (E) = 0$, by the item $(i)$. Thus $\pi_0 (F) = 0$, that is, the Milnor fiber $F$ is path connected.

\vspace{0.3cm}

\noindent $(iii)$ Consider the following piece of the long exact sequence in homotopy of the Milnor fibration (\ref{MasseyII}):
$$\xymatrix{
\pi_1(E) \ar[r]^{{(f_|)}_*} & \pi_1(S^1) \ar[r] & \pi_0(F) \ar[r] & \pi_0 (E) = 0.
}$$

Suppose $F$ is path connected, that is, $\pi_0(F) =0$. Then ${(f_|)}_*$ is surjective and, therefore, there is a map $\lambda: S^1 \to E$ such that $[id] = {(f_|)}_* ([\lambda]) = [f\circ \lambda]$, where $[id]$ is the homotopy class of the identity map. Thus $\lambda: S^1 \to E$ is a cross-section.

Conversely, suppose that the Milnor fibration admits a cross-section $\lambda: S^1 \to E$. By the definition $[id] = [f\circ \lambda] = {(f_|)}_* ([\lambda])$ and, therefore, ${(f_|)}_*$ is surjective, since $\pi_1(S^1)$ is the infinity cyclic group generated by $[id]$. Thus, since
$$\xymatrix{
\pi_1(E) \ar[r]^{{(f_|)}_*} & \pi_1(S^1) \ar[r] & \pi_0(F) \ar[r] & \pi_0 (E) = 0
}$$
is an exact sequence, we conclude that $\pi_0 (F) = 0.$

\cqd

\vspace{0.1cm}

The next example comes from M. Ribeiro's Ph.D thesis in \cite{mr}. It shows that for $p=2$ only the Milnor conditions $(a)$ and $(b)$ are not enough to guarantee the connectness of the Milnor tube fibers.

\begin{exam} Consider $f:(\mathbb{R}^{3},0)\to (\mathbb{R}^{2},0),$ $f(x,y,z)=(xy,xz).$ Note that $\Sigma_{f}=\{x=0\}$ and $V=\{x=0\}\cup \{y=z=0\},$ hence Milnor condition $(a)$ holds. Just in the same way one has $\Sigma_{(f,r)}=\{x=0\}\cup\{x^2=y^2+z^2\},$ $\overline{\Sigma_{(f,r)}\setminus V}=\{x^2=y^2+z^2\},$ and so Milnor condition $(b)$ also holds true. Then, Theorem \ref{mtu} assures that the projection

$$f_{|}:B_{\epsilon}^{3}\cap f^{-1}(B_{\eta}^{2}\setminus\{0\})\to B_{\eta}^{2}\setminus\{0\}$$
is a locally trivial fibration for all $0<\eta \ll \epsilon$ small enough. Now, since the space $\{x=0\}$ disconnect the closed ball $B_{\epsilon}^{3}$ into two halves one has that the Milnor tube $B_{\epsilon}^{3}\cap f^{-1}(S_{\eta}^{1})$ is disconnected, and so are the Milnor fibers.
\end{exam}

\section*{Acknowledgment}\label{ack}

The authors would like to express thanks to the anonymous referee for important comments which improved the presentation
of the paper. They would also like to thank Raimundo Ara\'ujo dos Santos for helpful discussions.

\section*{Conflict of Interest Statement}
The authors declare that there is no conflict of
interest.


\end{document}